%Paper: alg-geom/9405009
%From: Oliver Kuechle <kuechle@math.ucla.edu>
%Date: Thu, 19 May 1994 17:17:10 -0700 (PDT)

% Title:  Some properties of Fano manifolds that are zeros of sections in
%         homogenous  vector bundles over Grassmannians
%
% Author: Oliver K\"uchle
%
% etc.:  8 pages, AmS-TeX version 2.1

\documentstyle{amsppt}

\def\MP #1{{\text{I\! P}^#1}}

\def\RO #1{{\uppercase\expandafter{\romannumeral #1}}}
\def\ro #1{{\romannumeral #1}}
\def\PP {{\text{I\! P}}}

\def\sln {{\bold{SL}(n,{\Bbb C})}}
\def\lra {{\longrightarrow}}
\def\Lra {{\Longrightarrow}}
\def\upa {{\uparrow}}

\def\Wedge {\bigwedge}
\def\Wwedge {{\textstyle\Wedge}}

\def\phi {\varphi}
\def\ZZ {{\Bbb Z }}

\def\CC {{\Bbb C}}

\magnification1200
\baselineskip=15pt

\topmatter

\title
Some properties of  Fano manifolds that are zeroes of sections in
homogenous vector bundles over Grassmannians
\endtitle

\author
Oliver K\"uchle
\endauthor

\address
Department of Mathematics, UCLA, Los Angeles, CA 90024-1555, U.S.A.
\endaddress

\email
kuechle\@math.ucla.edu
\endemail



\rightheadtext{Properties of Fano manifolds}

\abstract
Let $X$ be a Fano manifold which is the zero scheme of a general global
section $s$ in an irreducible homogenous vector bundle over a Grassmannian.
We prove that the restriction of the Pl\"ucker embedding embeds $X$
projectively normal, and that every small deformation of $X$ comes
from a deformation of the section $s$.
These results are strengthened in the case of Fano 4-folds
\endabstract

\endtopmatter

\document

\head
INTRODUCTION
\endhead
Let $Gr(k,n)=\sln /P_k$ be the Grassmannian of $k$-dimensional quotients
of $n$-dimensional complex space $\CC^n$ considered as quotient of
$\sln$ by a maximal parabolic subgroup $P_k$. Then (irreducible)
representations of $P_k$ give rise to (irreducible) homogenous vector bundles
over $Gr(k,n)$. The purpose of this note is to prove  the following
theorems:
\proclaim {Theorem 1} Let $X$ be a Fano manifold which is the zero scheme of
a general global section in a globally generated irreducible homogenous
vector bundle ${\Cal F}$ over $Gr(k,n)$. Then $X$ is projectively normal.
\endproclaim
Here by a Fano manifold we mean a manifold $X$ with ample anticanonical
divisor $-K_X$, and $X\subset
Gr(k,n)$ is
considered to be embedded by the restriction of the Pl\"ucker embedding.
\proclaim {Theorem 2} Let $X$ be as
 above. Then every small deformation of
$X$ is again the zero scheme of a section in the same homogenous bundle.
\endproclaim
Moreover it becomes obvious from the proof that the bundle ${\Cal F}$ in
Theorem 1 can be replaced by the sum of one irreducible
vector bundle and line bundles.
\smallskip
Concerning Fano 4-folds we have a slightly more general result:
\proclaim {Theorem 3}
Suppose $dim( X)=4$ and that the Picard group $Pic(X)$ of $X$ is
 generated by ${\Cal O}_X(-K_X)$. Then the statements of Theorems 1 and 2
remain true is ${\Cal F}
$ is fully reducible, i.e. direct sum of irreducible homogenous vector bundles.
\endproclaim
The idea of the proofs goes back to Borcea \cite{Bor} and Wehler \cite{We}, who
obtained
results similar to Theorem 2 in the case of varieties parametrizing
linear subspaces on complete
intersections of hypersurfaces. Nevertheless it is difficult to make
statements for arbitrary irreducible homogenous vector bundles.
Therefore it is worthwhile to point out
that the condition of $X$ being Fano is crucial here and sharp in a certain
sense (cf. Example 4.11
and Remark 5.5).
On the other hand homogenous vector bundles seem to play an important role
in the classification of Fano manifolds (cf. \cite{Muk}, \cite{K\"u}), e.g.,
all
Fano 3-folds $V$ of the "main series", i.e. with very ample $-K_V$ and
$b_2(V)=1$ arise as sections of the sum of an irreducible homogenous vector
bundle and line bundles over ordinary or isotropic Grassmannians.
\smallskip
The proofs work as follows: The Theorems will follow from the vanishing of
certain cohomology groups of bundles and sheaves on $X$ and $Gr(k,n)$. Via
spectral sequence arguments this is reduced to the vanishing of cohomology
groups of vector bundles on $Gr(k,n)$ involving only wedge products
and tensor products of ${\Cal F}$, its dual and the tangent bundle
$\Theta_{Gr(k,n)}$. Since all these bundles are homogenous, we can
apply Bott's Theorem to obtain the vanishings once we determine the
weights of the corresponding representations. This last step is the
"ugly" part and consists in combining various estimates, and here is where
the Fano-condition comes in to keep control of the shape of the occuring
weights.
\medskip I would like to thank Rob Lazarsfeld for helpful discussions,
UCLA for hospitality and the Deutsche Forschungsgemeinschaft for
financial support.\bigskip\noindent
\head
{\bf 1.}
\endhead
\medskip\noindent
Let $Y=Gr(k,n)\subset\PP^N$ be embedded by the Pl\"ucker embedding.
It is well known that $Y$ is
projectively normal of dimension $k(n-k)$ with canonical line bundle
${\Cal O}_Y(K_Y)={\Cal O}_Y(-n)$.
Let $X\subset Y$ be a subvariety.
Then from the commutativity of
$$\matrix & & && 0 && 0 &&\cr
&&&&\upa &&\upa &&\cr 0 &\lra &{\Cal J}_X &
\lra &{\Cal O}_Y &\lra &{\Cal O}_X &
\lra &0\cr &&&& \upa && \upa &&
\cr &&&& {\Cal O}_{\PP^N} &= & {\Cal O}_{\PP^N} &&\cr
&&&&&& \upa &&\cr &&&&&& {\Cal I}_X &&\cr
&&&&&& \upa &&\cr &&&&&& 0 &&\cr\endmatrix
\leqno(1.1)$$
where ${\Cal J}_X$ is the ideal sheaf of $X$ in $Y$ and ${\Cal I}_X$ the
ideal sheaf of $X$ in $\PP^N$, it is clear that ${\Cal O}_X(1):={\Cal O}_Y
(1)\otimes {\Cal O}_X$ embeds $X$ projectively normal if and only if
$$H^1(Gr(k,n),{\Cal J}_X(r))=0 \quad\forall\quad r\ge 1.\leqno(1.2)$$
Now let $s\in H^0(Gr(k,n),{\Cal F})$ be a general global section
in a globally generated vector bundle ${\Cal F}$ over $Gr(k,n)$, such
that $X$, the variety of zeroes of $s$, is non-empty.
Then it is known (cf. \cite{We}) that every small deformation of $X$ can be
obtained by varying the section $s$ if
$$ H^1(Gr(k,n), {\Cal F}\otimes {\Cal J}_X)=0,
   \quad \text{and} \leqno(1.3.\text a)$$
$$ H^1(X, \Theta_{Gr(k,n)}|_X)=0. \leqno(1.3.\text b)$$
Finally the Koszul complex associated to the section $s$
gives, for any vector bundle ${\Cal E}$ on $Gr(k,n)$, spectral sequences
$$\matrix
H^p(Gr(k,n),{\Cal E}\otimes \Wedge^q{\Cal F}^*)
      &\Lra &H^{p-q}(X,{\Cal E}|X)&\cr &&&\cr
H^p(Gr(k,n),{\Cal E}\otimes \Wedge^{q+1}{\Cal F}^*) &\Lra &
H^{p-q}(Gr(k,n),{\Cal E}
\otimes {\Cal J}_X), & q\ge 0.\cr
\endmatrix
\leqno(1.4)$$
\bigskip
\head
{\bf 2.}
\endhead
\medskip
Now we recall the set-up and fix notations for the weight calculations. Let
$$P_k=\Biggl\{\left( \matrix
 A & 0 \cr B & C \cr
\endmatrix \right)\in \sln,
A\in \bold{GL}(k,\CC)\Biggr\},$$ such that $Gr(k,n)=\sln /P_k$.
Then an irreducible homogenous vector bundle ${\Cal F}$ comes from
a representation of the reductive part of $P_k$ consisting of matrices
with $B=0$. Such a representation is uniquely determined by its highest
weight which can be  written as an integral vector
$\beta=(\beta_1,...,\beta_n)$ (cf. \cite{K\"u}).
In this notation the universal quotient bundle ${\Cal Q}$ on the Grassmannian,
which comes from the representation of $A$, has highest weight $(1,0,...,0)$,
and the tangent bundle $\Theta_{Gr(k,n)}$ reads $(1,0,...,0,-1)$.
The vector bundle ${\Cal F}$ is globally generated if  and only if
$\beta_i\ge\beta_{i+1}$ for all $1\le i\le n-1$, and the highest
weight of the representation corresponding to the dual bundle ${\Cal F}^*$
is $(-\beta_k,...,-\beta_1,-\beta_n,...,-\beta_{k+1})$.\par
The Weil group acts on the weights by permutations among the first $k$ and
the last $n-k$ entries, yielding
$$rk({\Cal F})=dim(\beta)=
\prod_{1\le i<j\le
k}{j-i+\beta_i-\beta_j\over{j-i}}
\prod_{k+1\le s<t\le
n}{t-s +\beta_s -\beta_t\over{t-s}}\leqno(2.1)$$
for the rank of the corresponding vector bundle.
In this notation, the vanishing part of Bott's Theorem (cf. \cite{Bot}, Theorem
IV')
can be expressed as follows:
\proclaim {Theorem}
Let ${\Cal F}$ be an irreducible homogenous vector bundle over $Gr(k,n)$ with
highest weight $(\gamma_1,...,\gamma_n)$.
Then $H^p(Gr(k,n),{\Cal F})$ does not vanish if and only if all entries
of the vector
$(\alpha_1,...,\alpha_n)=(n+\gamma_1,n-1+\gamma_2,...,1+\gamma_n)$
are distinct and $p$ is the number of pairs $(i<j)$ such that
$\alpha_i<\alpha_j$.
\endproclaim
\bigskip\head
{\bf 3.}
\endhead
\medskip
Let $X$ and ${\Cal F}$ be as in the introduction. Since $X$ has positive
dimension, i.e. $rk({\Cal F})<k(n-k)$ and (2.1) we may assume that
${\Cal F}$ comes from a representation of the $A$-part of
$P_k$, which means $\beta_{k+1}=...=\beta_n=0$ for the corresponding highest
weight. Hence we will write $\beta=(\beta_1,...,\beta_k)$ for this weight, and
$|\beta|=\beta_1 +...+\beta_k$. Then by symmetry the weight of $det({\Cal F})$
is $rk(\beta)\cdot |\beta|/k$ times the weight $(1,...,1)$ which corresponds
to $\Wedge^k {\Cal Q}$ defining the Pl\"ucker embedding. Therefore, using
the adjunction formula, the condition of $X$ being Fano turns out
to be
$$ n>{{rk(\beta)\cdot |\beta|}\over{k}}\leqno(3.1)$$
which is usually stronger than $dim(X)>0$ which in turn reads
$$ n>{{rk(\beta)}\over{k}}+k.\leqno(3.2)$$
Note that (3.1) shows in particular that the entries of the weights of
arbitrary wedge
products of ${\Cal F}$ are smaller than $n$.
\medskip\noindent
\remark {(3.3) Remark}
Our results are well known for hypersurfaces in projective space, so we may
assume  $k,n-k\ge 2$. Moreover we assume  that $\beta$ is not
one of the following (cf. \cite{K\"u}):\smallskip
(i): $(1,0,...,0)$, since $X=Gr(k,n-1)$ in this case.
\smallskip
(ii): $(2,0,...,0)$, $2k\le n$, since $X$ parametrizes $k$-dimensional
subspaces
on (affine) $(n-1)$-dimensional quadrics in this case. Then $X$ is itself
a rational homogenous manifold (resp. two disjoint copies of those if $2k=n$),
hence known to be rigid and projectively normal.\smallskip

(iii): $(2,0,...,0)$ and $(1,1,0,...,0)$ for $2k>n$, since $X=\emptyset$ is
these cases.
\endremark
\bigskip
\head
{\bf 4.}
\endhead
\medskip\noindent
Now we prove Theorem 1. By (1.2) and (1.4) it suffices to show
$$H^p(Gr(k,n),({\textstyle\Wedge}^{p}{\Cal F}^*)(r))=0 \quad\forall \quad p>0,
r>0.$$
We will prove a slightly stronger result, namely the vanishing for
$p>0$ and $r\ge 0$, which also implies the connectivity of $X$.
Suppose one of these groups does not vanish and
let $\beta$ be the highest weight of ${\Cal F}$.
Note that $(\Wedge^p{\Cal F}^*)(r)$ is fully reducible and apply
Bott's theorem.  The proof relies on the elementary observation that, since
there is no space "between" the last $n-k$ entries, $p$ jumps in steps of
$n-k$ and there has to be a $(n-k)$-jump between the entries of a weight
of $\Wedge^p{\Cal F}$. The latter forces $\beta$ to have a jump which makes
the rank big. But the rank in turn is bounded in terms of $n$ which yields
a contradiction.
\par
More formally, there has to be a weight $(b_1,...,b_k)$ of $\Wedge^p{\Cal F}$
and a positive integer $s<k$ such that
$$\left\{ \matrix
 p=s(n-k)\cr
  \cr n-1\ge b_1\ge ...\ge b_s\ge n-k+r+s\cr
 \cr 0\le b_k\le ...\le b_{s+1}\le r+s\cr
\endmatrix\right. \leqno(4.1)$$
Moreover, starting with $s(n-k)|\beta|=b_1+...+b_k$, (4.1) together
with (3.1) implies a  condition on the rank of $\beta$:
$$ rk(\beta) |\beta |< {{k^2(k-1)}\over{s(|\beta|-1)}} +k^2 \leqno(4.2)$$
Beginning with some special cases we show that such bundles do not exist.
\par\noindent {\bf (4.3)}
Using (3.3) we   may assume that $\beta\ne (1,0,...,0),$ $ (2,0,...,0),$ $
(1,...,1,0)$
or $(t,...,t)$.
\par\noindent {\bf (4.4)}
Suppose $k=2$. Then $s=1$, $r=0$, and $b_2\le 1$, so $b_1+b_2\le n$, but
$p=n-2$, hence $\beta=(1,0),$ $(1,1)$, or $(2,0)$.\par
Suppose $k=3$. If $s=1$ then $r\le 1$, $b_2+b_3\le 4$. By $|\beta |(n-3)\le
n+3$, we know $| \beta |\le 4$, and $|\beta |\ge 3$ by (4.3). If $|\beta |=3$,
then $n\le 6$ and $\beta =(3,0,0)$ or $(2,1,0)$ contradicting $n>rk(\beta)$.
If $|\beta |=4$, then $n=5>4/3\cdot rk(\beta)$ shows $\beta =(2,1,1)$.\par
If $s=2$, then $r=0$, so $b_1+b_2+b_3=(2n-6)|\beta|\le 2n$, but $|\beta|\ge 3$,
hence $n\le 4$.\par
So we may assume $k\ge 4$.
\par\noindent {\bf (4.5)}
Using (4.4), we conclude $|\beta |\le k$, in particular $\beta_k=0$.
\par\noindent {\bf (4.6)}
We may assume $\beta\ne (c,0,...,0)$, since, for $c\ge 3$,
$$rk(c,0,...,0)={{k+c-1}\choose {c}}\ge
{{k^2(k-1)}\over{c(c-1)}}+{{k^2}\over{c}}.$$
Suppose $\beta =(1,1,0,...,0)$. By (3.3iii), $n-k\ge k$, but
$b_i\le k-1$ for the entries of wedge products of $\beta$.\par
Suppose $\beta =(2,1,0,...,0)$. Then $3\cdot rk(\beta)=k(k^2-1)\ge
k^2(k-1)/2+k^2$.
For $\beta =(2,2,0,...,0)$, we have
$4\cdot rk(\beta)=k^2(k^2-1)/3\ge k^2(k-1)/3+k^2$ since $k\ge 4$. Moreover
$rk(t,t,0,..,0)\ge
rk(2,2,0,...,0)$ for $t\ge 2$ and $rk(v,w,0,...,0)\ge rk(2,1,0,...,0)$ for
$v>w>0$. \par
Hence, by (4.2) we may assume $\beta_3 \ne 0$.
\par\noindent {\bf (4.7)}
In the same way it is shown that we may assume $\beta_{k-2}=0$ and $|\beta |\ge
5$.
\par\noindent {\bf (4.8)}
We are left with $k\ge 6$, $|\beta |\ge 5$  and
     $ rk(\beta)\ge {k\choose 3}$,
which gives a contradiction to (4.2).\smallskip
This completes the proof of Theorem 1.
\proclaim {(4.9) Corollary}
Under the assumptions of Theorem 1 $X$ is connected or $2k=n$ and
    ${\Cal F}\simeq S^2{\Cal Q}$.
\endproclaim
Considering a diagram similar to (1.1) and using Kodaira vanishing, we obtain
\proclaim {(4.10) Corollary}
For $X$ to be projectively normal it is sufficient to assume that ${\Cal F}$
is the sum of one irreducible homogenous vector bundle and line bundles.
\endproclaim

The following  example shows that one has to impose a non-positivity
condition on the canonical bundle.
\par\noindent
\example{(4.11) Example}
Consider the surface $S$ parametrizing lines on the cubic 3-fold, in this
context originally studied by Wehler \cite{We}. $S$ is the variety of zeroes
of a section in $S^3{\Cal Q}$ over $Gr(2,5)$.  Then one can show $H^3(Gr(2,5),
\Wedge^3(S^3 {\Cal F})^*(2))=1$ and $H^1(Gr(2,5), {\Cal J}_S(2))=1$,
i.e., $S$ is not quadratically normal with respect to ${\Cal O}_S(1)={\Cal
O}_S(K_S).$
\endexample
\bigskip
\head
{\bf 5.}
\endhead
\medskip
The proof of Theorem 2 is similar. By (1.3) and (1.4) it suffices to show
$$H^p(Gr(k,n),{\Cal F}\otimes \Wwedge^p{\Cal F}^*)=0\quad\forall \quad
p\ge  1\quad \text{and} \leqno(5.1.\text a)$$
$$H^{p+1}(Gr(k,n),\Theta_{Gr(k,n)}  \otimes \Wwedge^{p}{\Cal F}^*)
=0\quad\forall\quad p\ge 0. \leqno(5.1.\text b)$$
If (5.1a) does not hold, then there exists an integer $0<s<k$ and a weight
$(a_1,...,a_k)$ of ${\Cal F}^*\otimes \Wedge^{s(n-k)}{\Cal F}$ such that
$$\left \{ \matrix
  n-1\ge a_1\ge ...\ge a_s\ge n-k+s\cr\cr
  a_k\le ...\le a_{s+1}\le s.\cr
\endmatrix\right. \leqno(\text a)$$
Therefore we   get
$$n\le
{{1}\over{s(|\beta|-1)}}\bigl(s(k-s)+|\beta|(sk+1)-s+1\bigr),\leqno(5.2)$$
where $\beta$ as always is the highest weight of ${\Cal F}$. Concerning
violations of (5.1.b) the output of Bott's Theorem is more subtle since
the weight of $\Theta_{Gr(k,n)}$ is involved. Namely, if (5.1.b) does not hold,
then there exists an integer $0<s<k$ and a weight $(b_1,...,b_k)$ of
$\Wedge^p{\Cal F}$ such that either
$$\left\{ \matrix  p=s(n-k)-1\cr
  \cr n-1\ge b_1\ge ...\ge b_s\ge n-k+s+1\cr
  \cr b_k\le ...\le b_{s+2}\le s,\; b_{s+1}\le s+1,\quad\cr
\endmatrix\right. \leqno(\text b)$$
or
$$\left\{ \matrix p=s(n-k)-2\cr
  \cr n-1\ge b_1\ge ...\ge b_{s-1}\ge n-k+s\cr
 \cr n-k+s-1\le b_s\le n-k+s\cr
 \cr b_k\le ...\le b_{s+2}\le s,\; b_{s+1}\le s+1.\cr
\endmatrix\right. \leqno(\text{b'})$$
Note that $\beta$ is not of type $(1,1,0,...,0)$, since $b_1\le n-k-1$ in this
case (cf. (3.3iii)). So by (3.3) we may assume $|\beta |\ge 3$ in the
following.
Now (b) again implies (5.2), whereas from (b') we infer
$$\matrix n &\le
&{{1}\over{s(|\beta|-1)}}\bigl(s(k-s)+2-k+|\beta|sk+2|\beta|\bigr)\cr
&&\cr
&\le & \cases
k+2+{{k+3}\over{|\beta|-1}},\quad s=1\cr
             \cr 2k\quad\quad \text{for}\quad s\ge 2.
\endcases\cr
\endmatrix
 \leqno(5.3)$$
By considering increasing values of $|\beta|$ and using (3.1) we see that
the situation $s=1$ and $2k+1\le n \le k+2+{{k+3}\over{|\beta|-1}}$ can be
excluded.\smallskip
Hence, by (5.2) and (5.3) it remains to consider the case $n\le 2k$.
We do this by listing all possible highest weights $\beta$ and checking
case by case that (a), (b) and (b') can not be satisfied.
\proclaim {(5.4) Lemma}
Suppose $|\beta|\ge 3$ and
$${rk(\beta)\cdot |\beta|\over {k}} < n \le 2k.$$
Then $\beta=(\beta_1,...,\beta_k)$ is one of the following:
$$\matrix
          (\ro1):(t,...,t), t\le n-1&(\ro2): (2,1,...,1)\cr\cr
          (\ro{3}): (2,...,2,1), n=2k& (\ro4): (1,...,1,0) \cr\cr
          (\ro5): (1,1,1,0,0), 7\le n\le 10& (\ro6): (1,1,1,0,0,0), n=11,12\cr
\cr (\ro7): (1,1,1,1,0,0), n=11,12&\cr
\endmatrix $$
\endproclaim
\noindent
The Lemma's proof is obvious. \par
The exclusion of $(\ro1)$ is immediate and
$(\ro7)$ violates (5.2) and (5.3) since $s\ge 2$ in this case.
The remaining types are most economically dealt with by determining
the highest weights of the irreducible summands of the relevant
representations and comparing these to (a),(b) and (b').\par
This completes the proof of Theorem 2.
\medskip\noindent
\remark{(5.5) Remark} As an illustration of the well known fact that
K3 surfaces have nonalgebraic small deformations we have nonvanishing in
(5.1.b) for a quartic in $\MP3$, i.e. ${\Cal F}=S^4{\Cal Q}$ over $Gr(1,4)$.
\endremark
\bigskip
\head
{\bf 6.}
\endhead
\medskip\noindent
To prove Theorem 3 note that Fano 4-folds with $Pic(X)\simeq {\ZZ}\cdot K_X$
arising as zeroes of sections in fully reducible homogenous vector bundles
over Grassmannians have been classified in \cite{K\"u}. There are only few
cases that are not covered by Theorems 1 and 2, and verifying the
assertions in the above way poses no further problems.
\bigskip
\head
{\bf 7.}
\endhead
\medskip\noindent
Finally we ask some questions arising in this context.\smallskip
\noindent
\remark {(7.1) Question} Are there examples of Fano-manifolds
$V$ ($dim(V)\ge 4$)  with very ample $-K_V$ spanning $Pic(V)$ which are
not projectively normal ?
\endremark
\medskip\noindent
\remark{(7.2) Question} Do Theorems 1 and 2 hold  for
fully reducible ${\Cal F}$ ? \par
\endremark
The main problem here is that in general there occur representations of both
blocks of the reductive part of $P_k$ which makes the application of
Bott's Theorem difficult.
\bigskip\noindent
\head
References
\endhead
\medskip\noindent
\vskip 10pt\noindent
{\bf [Bor]} {\par\vskip-13pt\noindent\leftskip1.9cm Borcea, C.: {\it Deforming
varieties of $k$-planes of projective complete intersections.} Pacific J. of
Math. {\bf 143}, 25-36 (1990)\smallskip}\noindent
{\bf [Bot]} {\par\vskip-13pt\noindent\leftskip1.9cm Bott, R.: {\it Homogenous
vector bundles.}  Ann. of Math. {\bf 66}, 203-248 (1957)\smallskip}\noindent
{\bf [K\"u]} {\par\vskip-13pt\noindent\leftskip1.9cm K\"uchle, O.: {\it On Fano
4-folds and homogenous vector bundles over Grassmannians.}  to appear in Math.
Z.\smallskip}\noindent
{\bf [Muk]} {\par\vskip-13pt\noindent\leftskip1.9cm Mukai, S.: {\it Biregular
classification of Fano 3-folds and Fano manifolds of coindex 3.} Proc. Nat.
Acad. of  Sci. USA {\bf 83}, 3000-3003 (1989)\smallskip}\noindent
{\bf [We]} {\par\vskip-13pt\noindent\leftskip1.9cm Wehler, J.: {\it Deformation
of varieties defined by sections in homogenous vector bundles.} Math. Ann. {\bf
268}, 519-532 (1984)\smallskip}\noindent

\enddocument